\documentclass[12pt]{article}
\usepackage{latexsym}
\usepackage{amssymb}
 
\newcommand{\proof}{{\noindent \bf Proof. }}
\newtheorem{thm}{Theorem}

\newtheorem{lem}{Lemma}

\newcommand{\C}{{\cal C}}

\newcommand{\Q}{{\cal Q}}

\date{}

\begin{document}
\begin{titlepage}
\title{\bf INTERLOCKED PERMUTATIONS}
{\author{{\bf  G\'erard Cohen}
\thanks{deceased}
\\FRANCE
\and{\bf Emanuela Fachini}
\\{\tt fachini@di.uniroma1.it}
\\''La Sapienza'' University of Rome
\\ ITALY
\and{\bf J\'anos K\"orner}
\thanks{Department of Computer Science, University of Rome, La Sapienza, 
via Salaria 113, 00198 Rome, ITALY}
\\{\tt korner@di.uniroma1.it}
\\''La Sapienza'' University of Rome
\\ ITALY}} 

\maketitle
\begin{abstract}

The Shannon capacity of a graph with a countably infinite vertex set is easy to determine if  the graph's chromatic number and clique number are equal.  Even this problem becomes intriguing, however, if we ask for the asymptotic growth of the clique number of the subgraph induced on sequences of vertices without repetitions. A suitably normalized value of this clique number has a limit called permutation capacity. We consider related problems. Our results are purely combinatorial and in line with previous work about permutation capacity.

\end{abstract}
\end{titlepage}
\section{Permutations and graph capacity}

Malvenuto and K\"orner \cite{KM}, later joined by Simonyi  \cite{KMS} introduced a class of combinatorial problems about permutations representing a formal generalisation of Shannon's problem of the zero-error capacity of the discrete memoryless channel \cite{Sh}. As is well known, Shannon rephrased his famously unsolved information theory problem  
as a purely combinatorial question about finite simple graphs. The vertices of an arbitrary finite simple graph are interpreted as input symbols to a channel. The adjacency relation of vertices represents the distinguishability of symbols. Transmission of messages  occurs by the selection of  a sequence of symbols, and the common length of these sequences represents the duration of the transmission process.  Two sequences are distinguishable at the receiving end of the channel if they have a pair of distinguishable symbols somewhere in the same position. A channel code is  a set of pairwise distinguishable sequences. The aim is to determine the maximum size of such sets. At the core of Shannon's problem there is a more technical question about large sets of pairwise distinguishable sequences of length 
$n$ from a finite alphabet $\cal X$, all having the same {\it composition}. This means that all the sequences involved can be obtained from one another by permuting their 
coordinates. In other words, there is a probability distribution $P$ on $\cal X$ such that for each element $a \in {\cal X}$ the element $a$ occurs in each sequence 
$nP(a)$ times. In the case of infinite alphabets, a natural set of sequences "of the same composition'' is, for any integer $n$, the set of all the permutations of the set $[n]=\{1,2,\dots, n\}.$ 
Closest to the Shannon setup is the case when distinguishability is defined in terms of an infinite graph whose vertex set is the set of natural numbers \cite{KM}.  The common feature across the various generalisations 
of graph capacity (including those in \cite{GKV}) is that it is sufficient for distinguishability if a suitable condition holds in a single coordinate. If that condition is satisfied in some  coordinate  for a pair of sequences, then those sequences are distinguishable irrespective of what happens in other coordinates. From an information-theoretic viewpoint this means that the channel being modelled is {\em memoryless}. An attractive problem in this setup is the case of the graph whose vertex set is that of the natural numbers with  
edges between consecutive numbers. In \cite{KM} it is conjectured that the clique number of the graph induced on the permutations of $[n]$ in the aforementioned way is the middle binomial coefficient ${n \choose \lfloor n/2 \rfloor} .$

\medskip

The problems we are addressing in the present paper are of a slightly different nature; they do not correspond to memoryless channels. The motivation for these new problems comes from a conjecture of the third author, concerning what is called {\em reversing families of permutations} in \cite{DE}. Two permutations, $\sigma$ and $\tau$ of $[n]$ are called {\em reversing} if there exists a pair of distinct elements, $a\in [n]$ and $b \in [n]$ such that 
$$\sigma^{-1}(a)=\tau^{-1}(b)\quad \hbox{and} \quad \sigma^{-1}(b)=\tau^{-1}(a).$$
In other words, representing permutations as linear orders, the permutations $\sigma$ and $\tau$ are reversing if some pair $a$ and $b$ of elements of $[n]$ occupy in the two permutations the same pair of positions while their individual positions are exchanged. A set of permutations is called {\em reversing} if any two permutations from the set are reversing. 
As it is reported in Ellis \cite{DE}, K\"orner conjectured in 2010 that the size of a reversing family of permutations grows exponentially in $n.$ 
The best upper bound on the maximum size $F(n)$ of a reversing family of permutations of $[n]$, due to Cibulka \cite{C}, is
$$F(n)\leq  n^{n/2+O(n/\log_{2}n)}$$
The authors were not able to improve on this. Instead, we shall examine the rate of growth 
of a family of permutations satisfying somewhat similar conditions about the relation of the images of some pairs of elements in all pairs of permutations.  We will introduce three different {\em difference relations} for pairs of permutations. A relation is a difference relation if its validity implies that the corresponding permutations are different. Our 
relations are in terms of a pair of elements $a\in [n]$ and $b \in [n]$ whose four inverse images by the two permutations are all different. The relations are in terms of the natural order of these images. In case of a reversing pair we required a pairwise coincidence of the inverse images of $a$ and $b.$ Excluding such coincidences allows us to consider analogous but less binding relations in the hope of developing some intuition for the initial problem about reversing families. We will call these relations interlocking.

\section{Local disjointness}

We will say that the permutations $\sigma$ and $\tau$ are {\em locally disjoint} if there exist two distinct elements $a \in [n]$ and $b \in [n]$ such that
$$\sigma^{-1}(a) <  \sigma^{-1}(b) <  \tau^{-1}(a) < \tau^{-1}(b), $$
or 
$$\tau^{-1}(a) < \tau^{-1}(b) <  \sigma^{-1}(a) <  \sigma^{-1}(b).$$

We would like to determine the largest size $M(n)$ of a set of permutations of $[n]$ such that any pair of these permutations is locally disjoint. The key to our corresponding result Theorem \ref{thm:lodi} is the following

\begin{lem}\label{lem:gre}

Let $G(n)$ be the number of those permutations of $[n]$ that are not locally disjoint from the identity permutation. Then
$$G(n)\leq 5^n$$
\end{lem}

\proof

We will show that the permutations not locally disjoint from identity can be mapped injectively into quinary sequences of length $n,$ i. e., $n$-length sequences 
from an alphabet of 5 elements. Let $\sigma$ be a permutation of $[n]$ not locally disjoint from the identity. Let us consider a directed graph with vertex set $[n]$ in which there is a directed edge 
from $c$ to $d$ whenever $\sigma$ the function corresponding to the permutation $\sigma$ maps $c$ to $d.$ It is well--known that this directed graph is the vertex-disjoint union of (consecutively) oriented cycles. Let us call a directed edge 
$(x,y)$ {\it increasing} if $x<y$ and let us call it {\it decreasing} if $x>y.$ Likewise, we call a well-oriented path {\em increasing} ({\em decreasing}) if all its edges are increasing (decreasing). 
A {\em monotone path} is a path that is either increasing or decreasing. We  can partition the edges of a well-oriented cycle in the cycle decomposition of our permutation into maximal (non extendable) increasing and decreasing paths. Because of their maximality, they have to alternate in the decomposition. We call a cycle {\em simple} if its 
decomposition has only two monotone paths. 

We claim that all the directed cycles of $\sigma$ 
are simple. To see this, suppose to the contrary that $\sigma$ has a directed cycle that is not simple. 
Consider along the cycle all the maximal sets of consecutive edges pointing in the same direction. Without restricting generality we can suppose that there are more than one of them. (Were this not true, we would have at least two decreasing paths.) iAssign to any such increasing path the interval having the endpoints of the path as its endpoints.  
We claim that two such intervals cannot have a non-empty intersection unless they contain one another. If, to the contrary, we had two maximal increasing paths not containing 
each other and having a non-empty intersection, then these would contain {\em crossing} directed edges pointing in the same direction. We say that two directed edges pointing in the same direction are crossing  if one of them has an endpoint such that the two endpoints of the other edge are at its two different sides. 
Notice that two directed edges pointing in the same direction cannot cross, otherwise their 4 endpoints would define a configuration assuring that the identity and $\sigma$ are locally disjoint, a contradiction. 

It is clear that the cycle must contain at least one maximal increasing path which is not contained in any other maximal path.
Then, clearly, except for its extremal points, our increasing directed path is vertex-disjoint from the two paths starting from its extremal points and pointing in the opposite (decreasing) direction. Let $a<c$ be the endpoints of our increasing path.
Let us denote by $c'$ the other end-vertex of the uniquely defined directed edge starting in $c.$ Thus, since the cycle is supposed not to be simple, $c' \not=a,$  implying $c'<a.$ Likewise, let us denote by $a'$ the starting vertex of the decreasing edge pointing to $a.$ Then the edges $(c, c')$ and the edges $(a'a)$ are crossing, a contradiction.

Consider an arbitrary but fixed directed cycle of length at least 3 in the decomposition of $\sigma.$ Suppose that $[a,b]$ is the interval having the smallest and the largest vertex of the fixed cycle as its endpoints. We claim that if another directed cycle from the decomposition of $\sigma$ contains any vertex from the interval $[a,b]$, then it must be entirely contained in this interval. In fact, if this were not the case, then there would be a directed edge in the second cycle having an endpoint in the interval $[a,b]$ and another one outside of it. 
But then this edge would cross an edge from the fixed cycle in both directions. This amounts to two crossing edges going in the same direction, contradicting our hypothesis   that $\sigma$ is locally disjoint from the identity permutation.

The above allows us to encode a permutation with the structure we just described in the following manner. A vertex in $[n]$ is encoded by 1 if it is a fixed point in $\sigma.$
It is encoded by 2 if it is the smallest vertex in a cycle of the cyclic decomposition and by 3 if it is the largest element of a cycle. A vertex is denoted 4 if it belongs to the increasing path of a cycle and a 5 if it is on the decreasing path of a cycle. It is easy to check that such a quinary code of a permutation locally disjoint from the identity 
allows us to reconstruct its cyclic decomposition. In other words, the mapping that assigns a quinary code of length $n$ to these permutations is injective. 

\hfill$\Box$

It should be clear that the lemma remains true if the identity permutation is replaced in the statement by an arbitrary but fixed permutation.
As an immediate consequence of the lemma, we can now determine the asymptotics of $M(n).$ We have

\begin{thm}\label{thm:lodi}
$${n! \over 5^n}\leq M(n) \leq {n! \over {6^{\lfloor n/ 3\rfloor }}}$$

\end{thm}

\proof

The lower bound can be proved by a greedy algorithm. As a matter of fact, we can construct a family of pairwise locally disjoint permutations of $[n]$ by iteratively selecting any permutation to be in the family and then eliminating from the set of the remaining permutations all those which are not locally disjoint from the one we have just selected. 
By the lemma we just proved, this selection process will not stop before we have chosen at least as many permutations as the stated lower bound. 

To prove the upper bound, it is sufficient to realise that for $n=3$ no two permutations of $[n]$ are locally disjoint. Then for an arbitrary $n$ let us  divide $n$ into disjoint triples of consecutive numbers. Consider those permutations of $[n]$ in which all these triples remain fixed while within a triple the order of the elements is arbitrary. If $n$ is not divisible by 3, then the one or two largest numbers in $[n]$ will be left invariant by our permutations. The resulting at least $6^{\lfloor n/ 3\rfloor }$ permutations still satisfy the condition that no two of them are locally disjoint.  We can partition the set of all the permutations of $[n]$ into similar sets, grouping together the entries of the linear order representing the given permutation into disjoint and consecutive triples and permuting the 3 elements of each group in an arbitrary manner. In a family of  locally disjoint permutations no two permutations of the same class can be present. This gives the upper bound.

\hfill$\Box$

\section{Encapsulation}

We now define a very similar condition and show that a similar result holds also in this case. We say that the permutations $\sigma$ and $\tau$ of $[n]$ are 
{\em locally encapsulating} if there exist elements $a$ and $b$ in $[n]$ for which in both permutations $a$ precedes $b,$ and furthermore, the intervals of positions having as endpoints these $a$ and $b$ in the two permutations, i.e. the intervals $[\sigma^{-1}(a), \sigma^{-1}(b)]$ and  $[\tau^{-1}(a), \tau^{-1}(b)]$ strictly contain one another. Let us denote by $N(n)$ the maximum number of pairwise encapsulating permutations of 
$[n].$ We will show that $N(n)>(\alpha n)!$ for every $\alpha<1$ and $n$ sufficiently large. More precisely,

\begin{thm}\label{thm:enc}
For every $\alpha <1$ and sufficiently large $n$, we have
$$(\alpha n)!<N(n)\leq n! $$
\end{thm}

\proof

The upper bound is trivial. To show the lower bound, fix an arbitrary positive integer $k$ and for every multiple $n$ of $k$ partition $[n]$ into disjoint consecutive intervals of length ${n \over k}.$ Let $\Q$ be the set of all those permutations of $[n]$ which leave all the consecutive intervals of our partition invariant. 
We have 
\begin{equation}\label{eq:inv}
|\Q|=\left[\left({n\over k}\right)! \right]^k.
\end{equation}
By a greedy algorithm, we find a set $\C\subset \Q$ of permutations such that any pair of its member permutations differ in at least two different intervals of coordinates and furthermore,
\begin{equation}\label{eq:c}
|\C| \geq {\left[\left({n\over k}\right)! \right]^{k-1} \over k}.
\end{equation}
In fact, this is true since for any permutation from $\Q$ there are less than $k\left({n\over k}\right)!$ permutations differing from it in just one interval. Notice that 
$$ \left[\left({n\over k}\right)! \right]^{k-1} k^{-1}=$$
$$=\left[\left({n\over k}\right)! \right]^{-1} n! {{\left[\left({n\over k}\right)! \right]}^k \over n!}k^{-1}
=(n[1-k^{-1}])! \left [\left({n\over k}\right)! (n[1-k^{-1}])!\right]^{-1} n!  {[\left(n/k\right)!]^k \over n!} k^{-1}$$
$$\geq (n[1-k^{-1}])! {2^{nh({1 \over k})} \over n+1}k^{-(n+1)}\geq (n[1-k^{-1}])! {k^{-(n+1)} \over n+1},$$
where the first inequality follows from the well-known non-asymptotic estimate 
$${[\left(n/k\right)!]^k \over n!} \geq k^{-n}$$
alongside with the observation that the binomial coefficient ${n\choose x}$ is lower bounded by ${2^{nh(x/n)} \over n+1},$ where $h(t)=-t\log_{2}t-(1-t)\log_{2}(1-t)$ is the binary entropy function.
To conclude, we just have to verify that any pair of permutations from $\C$ is an encapsulating pair. To this end, notice that for any two permutations of the same ground set 
we have an element of the ground set that appears earlier in one of the permutations, and another one that appears later. Consider two arbitrary permutations from $\C,$
$\sigma$ and $\tau.$ Since $\sigma$ and $\tau$ differ on (at least) two intervals in $[n],$ they contain different permutations of the elements of these intervals. Therefore, we can find an element 
$a$ of the first interval which appears in $\sigma$ strictly earlier than in $\tau.$ Likewise, for the second interval we can find an element $b$ that comes in $\sigma$ strictly later 
than in $\tau.$ This implies that the interval of positions $[\sigma^{-1}(a), \sigma^{-1}(b)]$ having its endpoints in $a$ and $b$  in $\sigma$ strictly includes the corresponding interval in $\tau.$
\hfill$\Box$

\section{Parallelism}

A further difference relation of the same kind gives a similar result. We will say that the permutations $\sigma$ and $\tau$ of $[n]$ are locally parallel if there are two different elements $a$ and $b$ of $[n]$ such that in $\sigma$ the position of $a$ precedes that of $b$ while in $\tau$ it is $b$ that precedes $a$, and further the interval of positions 
of $a$ and $b$ in $\sigma$ is disjoint from the corresponding interval in $\tau.$ More precisely, 
$$[\sigma^{-1}(a), \sigma^{-1}(b)]\cap [\tau^{-1}(b), \tau^{-1}(a)]=\emptyset $$
We would like to determine the largest cardinality $P(n)$ of a set of pairwise locally parallel permutations of $[n].$ To this end, 
we will upper bound the number of permutations of $[n]$ that are not locally parallel to a fixed permutation.

\begin{lem}\label{lem:par}
Let $H(n)$ be the number of those permutations of $[n]$ which are not locally parallel to the identity permutation. Then
$$H(n) \leq 6^n.$$
\end{lem}

\proof

If a permutation $\sigma$ is locally parallel to identity, then there are two elements, $a<b$ in $[n]$ such that the interval of positions $[a,b]$ is disjoint from the interval of positions $[\sigma^{-1}(b), \sigma^{-1}(a)]$, 
where, also, $\sigma^{-1}(b)<\sigma^{-1}(a).$ This means that in the cyclic decomposition of $\sigma^{-1}$ the directed edge $(b, \sigma^{-1}(b))$ is running parallel with the directed edge $(a, \sigma^{-1}(a))$ 
without crossing it. It is therefore sufficient to prove that the number of those permutations $\tau=\sigma^{-1}$ for which no such parallel pairs of directed edges $(b, \tau(b))$, $(a, \tau(a))$ exist, is at most $6^n.$ Let us now consider a permutation $\sigma$ not locally parallel to the identity and let us look at its directed edges pointing in the direction of their larger endpoint. Let us encode with $x$ the initial vertices of these edges, by $y$ their endpoints and by $z$ the points that are both initial point and endpoint of an increasing edge. There is no pair of parallel edges among these if and only if for any pair of them the edge that begins earlier, ends earlier, too. This means that to any given 
sequence of positions of the $x$'s, the $y$'s and the $z$'s, there is only one configuration without parallel edges. The same consideration applies to the decreasing edges, 
with the symbols $x', y', z'$ to encode initial vertices, end vertices or initial and end vertices of the descending directed edges, respectively. This means that to any $n$-length sequence over the 6-element alphabet $\{x,y,z,x',y',z'\}$ there corresponds at most one permutation not locally parallel to identity.

\hfill$\Box$

This lemma immediately applies to prove

\begin{thm}\label{thm:six}
$$ {n! \over {6^n}}<P(n)\leq {n! \over {6^{\lfloor n/ 3\rfloor }}}$$
\end{thm}

\proof

The upper bound holds for the same reason as in Theorem \ref{thm:lodi}. The lower  bound can be proved by the greedy algorithm. We use Lemma \ref{lem:par} to 
notice that in each step of the greedy algorithm, we have to remove strictly less than $6^n$ permutations.

\hfill$\Box$

\section{Pattern avoidance}

Lemma \ref{lem:gre} and Lemma \ref{lem:par} are somewhat reminiscent, both in form and content, of the well-known problem area of pattern avoidance in permutations \cite{Bo}. We define a pattern 
on a fixed finite number of coordinates, and ask for the number of those permutations in which this pattern is absent. As in the classical pattern avoidance problems, we find 
that the number of those permutations that do not contain our pattern is only exponential. What we call a pattern is basically different from the usual concept of pattern in the 
sense of Stanley and Wilf \cite{Bo}. At this stage we prefer not to give a definition so as to remain vague and flexible.  It would be interesting to consider pattern avoidance in more generality so as to comprise our present problems. 

\section{Acknowledgement}
We are grateful to David Madore and Mikl\'os B\'ona for their interest in this work.

\end{document}